\documentclass[12pt]{article}%{book}a4paper,
\usepackage{amsfonts}
\begin{document}
\large

\begin{center}
{\Large{\bf On the Well-possedness  of the Problem of Reconstruction
of Non-separate Boundary Conditions}}
\end{center}

\vspace{0.2cm}

\centerline {Akhtymov A. M.$^{a}$, Mouftakhov A. V.$^{b}$, Teicher
M.$^{b}$, Yamilova L. S.$^{c}$ }

\vspace{0.2cm}

$^a$ Institute of mechanics, Ufa, Russia;

$^b$ Department of Mathematics, Bar-Ilan University, Ramat-Gan,
Israel;

$^c$ Department of Differential Equations,  Bashkir State
University, Ufa, Russia.

\vspace{0.2cm}

\date{}

\section{Introduction}
The identification of boundary conditions of a spectral problem is
an important practical problem.

The question arises whether one would be able to detect boundary
conditions, using finite number  eigenvalues. The following papers
give and substantiate a positive answer to this question for several
cases (see  \cite{Akhtyamov 03 IRAN}, \cite{Akhtyamov 03 AJ-a},
\cite{Akhtyamov 03 DU}, \cite{Akhtyamov 01}, \cite{Akhtyamov 03
IPE}). In this paper we continue these researches.

 The problem in question belongs to the class of inverse problems and
 is a completely natural problem of identification of the boundary
conditions.

The problem of determining a boundary condition has been considered
in \cite{Frikha 00}. However, as data for finding the boundary
conditions and as opposed to condensation and inversion (as in
\cite{Frikha 00}), we take a set of eigenvalues.

Similarly formulated problems also occur in the spectral theory of
differential operators, where it is required to determine the
coefficients of a differential equation and the boundary conditions
using a set of eigenvalues (for more details, see \cite{Borg 46},
\cite{Levitan 87}, \cite{ Marchenko 77}, \cite{Poshel 87},
\cite{Akhtyamov 99 DM}, \cite{Yurko 00}). However, as data for
finding the boundary conditions, we take one spectrum, not several
spectra or other additional spectral data (for example, the spectral
function, the Weyl function or the so-called weighting numbers),
that were used in these papers. Moreover, their principal aim was to
determine the coefficients in the equation, not in boundary
conditions. The aim of this work is to determine the boundary
conditions of the eigenvalue problem, from its spectrum in the case
of a known differential equation.

 We consider an inverse spectral problem with the
third-order differential equation and the non-separated boundary
conditions.

Two theorems on the uniqueness of the solution of this problem are
proved, and a method for establishing the unknown conditions is
obtained, using $19$ eigenvalues.

 The method
of approximate calculation of unknown boundary conditions is
explained, with the help of an example.

\section{Formulation of the Inverse Problem} The following spectral problem is
considered:
$$l (y) =y'''(x)  +   \lambda \, p_1\,  y'(x) +  \lambda^2 \,p_2\,
y'(x) + \lambda^3  \, p_3 \, y(x) = 0, \eqno{(1)}$$
$$U_i(y)=\sum_{k=1}^3 \left( a_{ik}\, y^{(k-1)}(0)+a_{i\,k+3}\,
y^{(k-1)}(1)\right) =0 \qquad (i=1,2,3),
 \eqno{(2)}$$
\noindent where $ a_{ik}$ and $p_{i}$, are not dependent on
parameter $\lambda $, and  $ a_{ik}, \, p_{i} \in \mathbb{C}$.

We denote the matrix formed by the coefficients $a_{ij}$, by $A$,
and its third-order minors, by $M_{ijk}$:
$$ A=\left\|
\begin{array}{ccccc}
a_{11} & a_{12} & a_{13} & \dots & a_{16}
\\
a_{21} & a_{22} & a_{23} & \dots & a_{26}\\
a_{31} & a_{32} & a_{33} & \dots & a_{36}\\
\end{array}\right\|,
\qquad M_{ijk}=\left|
\begin{array}{ccc}
a_{1i} & a_{1j} & a_{1k}\\
a_{2i} & a_{2j} & a_{2k} \\
a_{3i} & a_{3j} & a_{3k} \\
\end{array}\right|$$
$\hspace{1.5 cm}(1\leq i<j< k\leq 6).$

In terms of problem (1)--(2), the inverse problem of reconstructing
the boundary conditions (2) can be stated as follows:

\vspace{0.1cm}

{\it Inverse problem.} The coefficients $a_{ij}$ of the forms $
U_i(y),\; (i=1,2,3)$ in problem (1)--(2), are unknown. The rank of
the matrix $A$ formed by these coefficients is equal to $3$. The
non-zero eigenvalues $\lambda_{k}$ of problem~(1)--(2), are known.
It is required to find the boundary condition, i.e. to reconstruct
the linear span $\langle \bf{a_{1}},\bf{a_{2}},\bf{a_{3}}\rangle$ of
vectors
${\bf{a_{i}}}=\left(a_{i1},a_{i2},a_{i3},a_{i4},a_{i5},a_{i6}\right),\,(i=1,2,3)$.

\section{Uniqueness Theorem for the Solution of the Inverse Problem}

Along with the problem (1)--(2), consider the following spectral
problem
$$l (y) =y'''(x)  +   \lambda \, p_1\,  y'(x) +  \lambda^2 \,p_2\,
y'(x) + \lambda^3  \, p_3 \, y(x) = 0, \eqno{(3)}$$
$$
\widetilde{U}_i(y)=\sum_{k=1}^3 \left( b_{ik}\,
y^{(k-1)}(0)+b_{i\,k+3}\, y^{(k-1)}(1)\right) =0 \qquad (i=1,2,3).
 \eqno{(4)}$$
We denote the matrix formed by the coefficients $b_{ij}$, by $B$,
and its third-order minors, by $\widetilde{M}_{ijk}$:
$$B=\left\|
\begin{array}{ccccc}
b_{11} & b_{12} & b_{13} & \dots & b_{16}
\\
b_{21} & b_{22} & b_{23} & \dots & b_{26}\\
b_{31} & b_{32} & b_{33} & \dots & b_{36}\\
\end{array}\right\|,
\qquad \widetilde{M}_{ijk}=\left|
\begin{array}{ccc}
b_{1i} & b_{1j} & b_{1k}\\
b_{2i} & b_{2j} & b_{2k} \\
b_{3i} & b_{3j} & b_{3k} \\
\end{array}\right|.$$

{\bf Theorem~ 1} (on uniqueness of the solution of the inverse
problem). {\it Suppose that the following conditions are satisfied:

$ {\rm rank}\, A= {\rm rank}\, B =3, $ not every number $ \lambda $,
is an eigenvalue for problems

(1)--(2), (3)--(4), the eigenvalues of the problem (1)--(2), and the
eigenvalues of the problem (3)--(4), coincide with their
multiplicities taken into account, and roots of the characteristic
equation $\omega^3 + p_1\,\omega^2 + p_2\,\omega + p_3 = 0$, satisfy
conditions:

1) $\sum_{i\in I_k}e_i\omega_i\neq 0 $, where $e_i=\pm1$ and $I_k$
is any subset of the set $\{1,2,3\},$

2) $ p_2\neq 0 $.

Then $Span\langle \bf{a_{1}},\bf{a_{2}},\bf{a_{3}}\rangle$ $\,=\,$
$Span\langle \bf{b_{1}},\bf{b_{2}},\bf{b_{3}}\rangle$, where
${\bf{a_{i}}}=\left(a_{i1},a_{i2},a_{i3},a_{i4},a_{i5},a_{i6}\right),$
${\bf{b_{i}}}=\left(b_{i1},b_{i2},b_{i3},b_{i4},b_{i5},b_{i6}\right),$
$\,(i=1,2,3)$.}

 {\bf Remark.} {\it It follows from the theorem's conditions that $p_1 \neq 0, \, p_3 \neq 0$}.

 {\bf Proof.}

By definition, put  $
z_{j}=\left\{\begin{array}{c}y^{(j-1)}(0),\;j=1,2,3;\\y^{(j-4)}(1),\;j=4,5,6.\end{array}\right.
$\\
%Тогда задачи (1)--(2), (3)--(4) можно переписать в виде
Then problems (1)--(2), (3)--(4), can be represented as
$$ l(y) = 0,\qquad U_i (y) = \sum_ {j=1} ^6 a_{ij} \, z_j = 0 \quad(i = 1, \, 2, \,
3). \eqno {(5)}$$
$$ l(y) = 0,\qquad \widetilde{U}_i (y) = \sum_ {j=1} ^6 \widetilde{a}_{ij} \, z_j =
0 \quad(i = 1, \, 2, \, 3). \eqno {(6)}$$

Let $\{ y_k(x)\}_{k=1,2,3}$, be a fundamental system of the solution
of equation (1), which meet the following conditions:

$y_k^{j-1}(0)=\left\{\begin{array}{c}1,\;k=j,\\0,\;k\neq
j\end{array}\right. \qquad (k,j=1,2,3). $

By definition, put $z_{kj}=\left\{\begin{array}{c}y_k
^{(j-1)}(0),\;j=1,2,3,\\y_k ^{(j-4)}(1),\;j=4,5,6\end{array}\right.
\qquad (k=1,2,3). $

By $\Delta(\lambda)$, denote the characteristic determinant of
problem (5). By $\widetilde{\Delta}(\lambda)$, denote the
characteristic determinant of problem (6).

 Let's consider the following boundary problem
$$ l(y) = 0,\qquad \widetilde{U}_i (y) = \sum_ {j=1} ^6 b_{ij} \, z_j = 0 \quad(i = 1, \,
2), \eqno {(7)}$$
$$ e^{f(\lambda)} \,\widetilde{U}_3 (y)\, = \, e^{f(\lambda)}\sum_ {j=1} ^6 b_{3j} \,
z_j, \eqno {(8)}$$ where $f(\lambda)$ is an entire function, which
will be chosen later.

The following function is a characteristic determinant of the
problem (7)--(8).

$$\widetilde{\Delta_1}(\lambda)=\left| \begin {array}{ccc}
\widetilde{U}_1(z_1)& \widetilde{U}_1(z_2)& \widetilde{U}_1(z_3)\\
\widetilde{U}_2(z_1)& \widetilde{U}_2(z_2)& \widetilde{U}_2(z_3)\\
e^{f(\lambda)}\widetilde{U}_3(z_1)& e^{f(\lambda)}\widetilde{U}_3(z_2)& e^{f(\lambda)} \widetilde{U}_3(z_3)\\
\end {array}\right|\equiv e^{f(\lambda)}\widetilde{\Delta}(\lambda).\eqno {(9)}$$
It immediately follows that the eigenvalues of spectral problems (6)
and (7) -- (8), are equal.

Hence, by the condition of the theorem, the eigenvalues of spectral
problems (5) and (7) -- (8), are equal as well.

$\Delta(\lambda)$ and $\widetilde{\Delta}_1(\lambda)$ are entire
functions of $\lambda $ (see \cite{Naimark}). The eigenvalues of
spectral problems (5), are zeros of function $\Delta(\lambda)$, and
the eigenvalues of spectral problems (7) -- (8) are zeros of
function $\widetilde{\Delta}_1(\lambda)$.
 Since the eigenvalues of
spectral problems (5) and (7) -- (8) are equal, zeros of functions
$\Delta(\lambda)$ and $\widetilde{\Delta}_1(\lambda)$, are equal
too. Hence, by consequence of Weierstrass' factorization theorem,
about zeros of an entire function \cite{Gol'dberg 91},
$$\widetilde{\Delta}_1(\lambda)\equiv C e\,^{g(\lambda)}\Delta(\lambda), \eqno {(10)}$$
where $g(\lambda)$ is an entire function.

It follows from (10) that
$$ e\,^ {f(\lambda)}
\left| \begin {array}{ccc}
U_1(y_1)& U_1(y_2)& U_1(y_3)\\
U_2(y_1)& U_2(y_2)& U_2(y_3)\\
U_3(y_1)& U_3(y_2)& U_3(y_3)\\
\end {array}\right|\equiv C e\,^ {g(\lambda)}
\left| \begin {array}{ccc}
\widetilde{U}_1(y_1)& \widetilde{U}_1(y_2)& \widetilde{U}_1(y_3)\\
\widetilde{U}_2(y_1)& \widetilde{U}_2(y_2)& \widetilde{U}_2(y_3)\\
\widetilde{U}_3(y_1)& \widetilde{U}_3(y_2)& \widetilde{U}_3(y_3)\\
\end {array}\right|.\eqno {(11)}$$

$f(\lambda)$ is an entire function. Let's choose $f(\lambda)$, so
that $f(\lambda)\equiv g(\lambda).$ Then it follows from (11) that
$$\left| \begin {array}{ccc}
U_1(y_1)& U_1(y_2)& U_1(y_3)\\
U_2(y_1)& U_2(y_2)& U_2(y_3)\\
U_3(y_1)& U_3(y_2)& U_3(y_3)\\
\end {array}\right|\equiv
C \left| \begin {array}{ccc}
\widetilde{U}_1(y_1)& \widetilde{U}_1(y_2)& \widetilde{U}_1(y_3)\\
\widetilde{U}_2(y_1)& \widetilde{U}_2(y_2)& \widetilde{U}_2(y_3)\\
\widetilde{U}_3(y_1)& \widetilde{U}_3(y_2)& \widetilde{U}_3(y_3)\\
\end {array}\right|.\eqno {(12)}$$

Taking into consideration the boundary conditions of the problems
(5) and (6), It follows from (12) that
$$\sum_ {i=1}^6\sum_ {j=1}^6\sum_ {k=1}^6\,z_{1i}z_{2j}z_{3k} M_{i\,j\,k}=
C \sum_ {i=1}^6\sum_ {j=1}^6\sum_ {k=1}^6\,z_{1i}z_{2j}z_{3k}
\widetilde{M}_{i\,j\,k}.\eqno{(13)}$$

It follows from (13) that
$$\sum_ {1\leq i<j<k\leq 6}\,\left(M_{i\,j\,k}-C \widetilde{M}_{i\,j\,k}\right)\,Z_{i\,j\,k}=0, \,\eqno{(14)}$$
where $Z_{i\,j\,k}= \left|
\begin{array}{ccc}
z_{1i} & z_{1j} & z_{1k}\\
z_{2i} & z_{2j} & z_{2k}\\
z_{3i} & z_{3j} & z_{3k}\\
\end{array}\right|.$

$\{Z_{i\,j\,k}|1\leq i<j<k\leq 6\}$ is a system of linearly
independent functions of $\lambda $. Then it follows from (14) that
$(M_{i\,j\,k}-C \widetilde{M}_{i\,j\,k})=0$, i.e. $M_{i\,j\,k}=C
\widetilde{M}_{i\,j\,k}$. Then $Span\langle
\bf{a_{1}},\bf{a_{2}},\bf{a_{3}}\rangle$ $\,=\,$ $Span\langle
\bf{b_{1}},\bf{b_{2}},\bf{b_{3}}\rangle$ (see \cite{Hodge 94}). This
completes the proof of the theorem.

\section{Exact Solution of the Inverse Problem}

Let us assume that $\lambda_k,\;(k=1,\ldots,19)$, are eigenvalues of
spectral problems (1) -- (2)

$$
y(x,\lambda)=C_1y_1(x,\lambda)+C_2y_2(x,\lambda)+C_3y_3(x,\lambda),\qquad
\eqno (15)
$$
where $\{y_i(x)|i=1,2,3\}$ is a fundamental system of the solution
of equation (1)--(2).

We shall find $C_i,\;(i=1,2,3)$, using the boundary condition (2).
We substitute (15) with (2). From this, we obtain the following
system of equations:
$$\begin{array}{c}
C_1U_1(y_1) + C_2U_1(y_2) + C_3U_1(y_3)=0,\\
C_1U_2(y_1) + C_2U_2(y_2) + C_3U_2(y_3)=0, \\
C_1U_3(y_1) + C_2 U_3(y_2)+ C_3 U_3(y_3)=0. \\
\end{array}$$
The non-zero solution for $C_i $ exists if and only if the following
determinant

$$\Delta(\lambda)=\left|\begin{array}{ccc}
U_1(y_1) & U_1(y_2) & U_1(y_3)\\
U_2(y_1) & U_2(y_2) & U_2(y_3) \\
U_3(y_1) & U_3(y_2) & U_3(y_3) \\
\end{array}\right|\qquad\eqno(16)$$
 is equal to zero \cite{Naimark}.

Transforming (16) and using the designations entered at the proof of
Theorem 1, we obtain
$$
\Delta(\lambda)=\sum_{1\leq i<j<k\leq 6}
Z_{ijk}(\lambda)M_{ijk}=0,\, \mbox{where}\, Z_{ijk}=\left|
\begin{array}{ccc}
z_{1i} & z_{1j} & z_{1k}\\
z_{2i} & z_{2j} & z_{2k} \\
z_{3i} & z_{3j} & z_{3k} \\
\end{array}\right|.\eqno(17)$$

Let's substitute values $ \lambda_m \; (m=1, \ldots, 19) $ to $
\Delta (\lambda) $, we shall receive the system with 19 homogeneous
equations of 20 variables $M_{ijk}$.
$$
\sum_{1\leq i<j<k\leq 6}
Z_{ijk}(\lambda_m)M_{ijk}=0.\qquad\eqno(18)$$

The system (18) has infinite number of solutions. If the rank of the
system equals 19, then  unknown minors $M _ {ijk} $ are determined
from the system, accurate to coefficient. By these minors, the
corresponding boundary conditions can be unequivocally found by
means of known methods of linear algebra.

Indeed, since $rank A = 3$ then one of minors $M _ {ijk} $ is not
equal to zero. Let $M _ {135} \neq0, $ then
%it is possible for
%matrix $A$ to lead, by means of linear transformations, to the
%following form:
after linear transformations the matrix $A$ can be written as
follows:

$$A=\left\|
\begin{array}{cccccc}
1 & a_{12} &0 & a_{14} & 0 & a_{16}
\\
0 & a_{22} &1 & a_{24} & 0 & a_{26}\\
0 & a_{32} & 0 & a_{34} & 1 & a_{36}\\
\end{array}\right\|.$$

With that, the minors $M _ {ijk} $ will not exchange (or probably
will be multiplied by non-zero number).

For this matrix, we get
 $$M_{135}=1,\,M_{134}=a_{34},\,M_{136}=a_{36},\,M_{123}=-a_{32},\,M_{235}=a_{12},$$
$$M_{156}=-a_{26},\,M_{145}=a_{24},\,M_{345}=-a_{14},\,M_{125}=a_{22},\,M_{356}=a_{16}.$$
Then matrix $A$ can be written as follows:
$$A=\left\|
\begin{array}{cccccc}
M_{135} & M_{235} &0 & -M_{345} & 0 & M_{356}
\\
0 & M_{125}/M_{135} &1 & M_{145}/M_{135} & 0 & -M_{156}/M_{135}\\
0 & -M_{123}/M_{135} & 0 & M_{134}/M_{135} & 1 &
M_{136}/M_{135}\\\end{array}\right\|\qquad\eqno(19)$$

This reasoning proves:

{\bf Theorem 2 }(on the uniqueness of the solution of the inverse
problem). {\it If the matrix of system (18) has a rank of $19$, the
solution of the inverse problem of the reconstruction boundary
conditions (2) is unique.}

\section{Example}
We shall consider application of a method of definition of the
boundary conditions by $19$ eigenvalues for the following boundary
problem:
$$ l(y)=y'''(x)-(3i+3)\lambda y''(x) +(9i-2)\lambda^2
y'(x)+6\lambda^3y(x)=0,\qquad\eqno (20)$$
$$U_i(y)=\sum_{k=1}^3 \left( a_{ik}\, y^{(k-1)}(0)+a_{i\,k+3}\,
y^{(k-1)}(1)\right) =0,\,(i=1,2,3).\qquad\eqno (21)$$

Let us know 19 eigenvalues of a problem (20) - (21)
$$ \lambda_1=0.46 - 0.12\times
i,\; \lambda_2=5.88+ 3.86\times i,\;\lambda_3=6.51 - 0.55 \times
i,$$ $$ \lambda_4=12.81 - 0.56 \times i,\;\lambda_5=19.1 - 0.56
\times i,\; \lambda_6=-4.27 + 0.51 \times i,$$
$$\lambda_7=-7.16 + 1.06\times
i,\; \lambda_8=-10.54 + i,\;\lambda_9= -13.50 + 1.32 \times i,$$
$$ \lambda_{10}=-19.81 + 1.49
\times i,\; \lambda_{11}=-23.1 + 1.41\times i,\;\lambda_{12}=-26.11
+ 1.61\times i,$$
$$\lambda_{13}=-29.38 + 1.54 \times
i,\;\lambda_{14}=-32.41 + 1.71 \times i,\;\lambda_{15}= -35.67 +
1.64 \times i,$$ $$ \lambda_{16}= -38.7 + 1.8 \times
i,\;\lambda_{17}=-44.99 + 1.87 \times i,\; \lambda_{18}=-48.23 +
1.80 \times i,$$
$$\lambda_{19}=-51.28 + 1.93 \times
i.$$

The fundamental system of decisions of the equation (20) satisfying
conditions
$$y_k^{j-1}(0)=\left\{\begin{array}{c}1,\;k=j\\0,\;k\neq j
,\end{array}\right. \qquad (k,j=1,2,3)$$ has the following
appearance:

$y_1=C_1e^{i\lambda x}+C_2e^{2i\lambda x}+C_3e^{3\lambda x},\quad$
$y_2=K_1e^{i\lambda x}+K_2e^{2i\lambda x}+K_3e^{3\lambda x},\quad$

$y_3=N_1e^{i\lambda x}+N_2e^{2i\lambda x}+N_3e^{3\lambda x},$

where
$$ C_1=\frac{9}{5}+\frac{3}{5}\,i,\quad
C_2=-\frac{9}{13}-\frac{6}{13}\,i,\quad
C_3=-\frac{7}{65}-\frac{9}{65}\,i,$$
$$
K_1=(-\frac{9}{10}+\frac{7}{10}\,i)/ \lambda,\quad
K_2=(\frac{9}{13}-\frac{7}{13}\,i)/\lambda,\quad
K_3=(\frac{27}{130}-\frac{21}{130}\,i)/\lambda,$$
$$
N_1=(\frac{1}{10}-\frac{3}{10}\,i)/ \lambda^2,\quad
N_2=(-\frac{2}{13}+\frac{3}{13}\,i)/\lambda^2,\quad
N_3=(\frac{7}{130}+\frac{9}{130}\,i)/\lambda^2.$$

Having solved system (18) by Maple, we shall find minors $M_{ijk}$:

$M_{135} = C,\;$ $M_{236} = (-5.65 - 3.95\,i)\times 10^{-19}C,\;$

$M_{134} = 0.50\,C -2.00\,i\times 10^{-8}C,\;$ $M_{145} =
(5.05+3.62\,i)\times10^{-8}C,\; $

$M_{136} = (-5.04-3.30\,i)\times 10^{-10}C,\;$ $M_{126} =
(-7.73+7.10\,i)\times10^{-10}C,\;$

$M_{156} = (2.44+ 0.56\,i)\times 10^{-7}C,\;$ $M_{256} =
(-0.03-1.21\,i)\times 10^{-8}C,\;$

$M_{356}= C+1.10\,i\times 10^{-8}C,\;$ $M_{346} =
0.5\,C-6.86\,i\times10{-9}C,\;$

$M_{124} = (1.24-5.13\,i)\times10^{-8}C,\;$ $M_{456} =
(-1.84-0.52\,i)\times10^{-7}C,\;$

$M_{125} = (- 2.17-3.30\,i)\times 10^{-8}C,\;$ $M_{146} = (0.85 -
3.79\,i)\times 10^{-8}C,\;$

$M_{235} = C +1.13\,i\times10^{-8}C,\;$ $M_{123} = (1.29+
1.33\,i)\times 10^{-9}C,\;$

$M_{345} = -0.5\,C -8.10\,i\times10^{-8}C,\;$ $M_{234} =
0.5\,C-8.20\,i\times10^{-9}C,\;$

$M_{245} = (0.63- 2.74\,i)\times10^{-8}C,\;$ $M_{246} = (6.63
-1.27\,i)\times10^{-8}C.$

Let $C=1$, then we have

$M_{125}\approx0,\;M_{145}\approx0,\;M_{156}\approx0,\;M_{123}\approx0,\;M_{136}\approx0$.

By (19), we get
$$A=\left\|
\begin{array}{cccccc}
1 & 1 &0 & 0.5& 0 &1
\\
0 & 0&1 & 0& 0 & 0\\
0 & 0 & 0 & 0.5 & 1 &0\\\end{array}\right\|.$$
Thus, unknown
boundary conditions are found.
$$U_1(y)=y(0)+y'(0) + 0.5\,y(1)+ y''(1)=0,$$
$$U_2(y)=y''(0)=0,\quad U_1(y)=0.5\,y(1)+ y'(1)=0.$$

\section{Acknowledgements}

This research was partially supported by the Russian Foundation for
Basic Research (06-01-00354a),  Emmy Noether Research Institute for
Mathematics, the Minerva Foundation of Germany, the Excellency
Center "Group Theoretic Methods in the Study of Algebraic Varieties"
 of the Israel Science Foundation, and by EAGER
(European Network in Algebraic Geometry).

\end{document}